\documentclass[12pt]{article}

\usepackage{amssymb}

\usepackage{pictex}
\usepackage{amsmath}
\usepackage{amscd}
\usepackage{amstext}
\usepackage{amssymb}
\usepackage{amsfonts}
\usepackage{latexsym}
\usepackage[dvips]{color}
\usepackage{tikz-cd}
\newtheorem{T}{Theorem}
\newtheorem{Lem}{Lemma}
\newtheorem{Prop}{Proposition}
\newtheorem{Def}{Definition}
\newtheorem{Cor}{Corollary}
\newtheorem{Rem}{Remark}

\newcommand{\bt}{\begin{T}}
	\newcommand{\bl}{\begin{Lem}}
		\newcommand{\bp}{\begin{Prop}}
			\newcommand{\bc}{\begin{Cor}}
				\newcommand{\bd}{\begin{Def}}
					\newcommand{\br}[2]{\begin{Rem}\label{#1}{\rm #2}}
						\newcommand{\er}{ \end{Rem}}
					\newcommand{\et}{\end{T}}
				\newcommand{\el}{\end{Lem}}
			\newcommand{\ep}{\end{Prop}}
		\newcommand{\ec}{\end{Cor}}
	\newcommand{\ed}{\end{Def}}
\DeclareMathOperator{\rank}{rank}
\newcommand{\be}{\begin{equation}}
	\newcommand{\ee}{\end{equation}}
\newcommand{\beq}{\begin{eqnarray}}
	\newcommand{\eeq}{\end{eqnarray}}
\newcommand{\beqq}{\begin{eqnarray*}}
	\newcommand{\eeqq}{\end{eqnarray*}}
\def\supp{\mathop{\rm supp \,}\nolimits}

\def\dist{\mathop{\rm dist \,}\nolimits}

\def\id{\mathop{\rm id \,}\nolimits}

\def\Id{\mathop{\rm Id}\nolimits}

\newcommand{\bpr}{{\bf Proof.}\ }
\newcommand{\epr}{\hspace*{\fill}\rule{3mm}{3mm}\\}
\newcommand{\cn}{{\mathcal N}}
\newcommand{\bspq}{B^s_{p,q}}

\newcommand{\Rd}{{\mathbb R}^{d}}

\newcommand{\R}{{\mathbb R}}
\newcommand{\N}{{\mathbb N}}
\newcommand{\Z}{\mathbb Z}

\newcommand{\C}{{\mathbb C}}

\newcommand{\cl}{{\mathcal{L}}}

\begin{document}

	\title{Embeddings  of  block-radial functions -- approximation properties and nuclearity }
	\author{ Alicja Dota\thanks{This research was supported  by the Pozna\'n University of Technology,  Grant no.  0213/SBAD/0118}, Leszek Skrzypczak} 
	\date{ }

	\maketitle


	\begin{abstract}
		 Let $ R_\gamma B^{s}_{p,q}(\Rd)$ be a subspace of the Besov space $B^{s}_{p,q}(\Rd)$ that consists of block-radial (multi-radial)  functions. We study an asymptotic behaviour of approximation numbers of compact embeddings $\id: R_\gamma B^{s_1}_{p_1,q_1}(\R^d) \rightarrow  R_\gamma B^{s_2}_{p_2,q_2}(\R^d)$.  Moreover we find the sufficient and necessary condition for nuclearity of the above embeddings. Analogous results are proved for fractional Sobolev spaces  $R_\gamma H^s_p(\Rd)$. 
		
		 
		 MSC[2020] 41A46, 46E35,  47B10 
	\end{abstract}
	


\section{Introduction}
 It is well known that symmetry  can generate compactness of Sobolev type embeddings on $\R^d$.    This  was noticed  in the case of the  first order Sobolev spaces of radial functions by W.~Strauss in  the seventies of the last century, cf. \cite{Strauss}. He also  proved  a pointwise estimate for  Sobolev radial functions that is call  Strauss' inequality now.  To the best of our knowledge  more general block-radial (multiradial) symmetry was first considered by L.P.~Lions in \cite{Lio82}. Compactness of embeddings in  both cases, radial and multi-radial, find applications in PDE's cf. e.g.   \cite{APS,Ku,PK,Serra}. 
  Later  many authors studied subspaces consisted of radial elements in different function spaces cf. \cite{EF1, SS1, SSV, SYY}. In particular asymptotic behaviour of entropy numbers of compact Sobolev embeddings of radial subspaces of Besov and Triebel-Lizorkin spaces were described  by {Th.~K{\"u}hn}, {H.-G.~Leopold, W.~Sickel} and the second named author in \cite{KLSS}, corresponding approximation numbers  was studied in \cite{ST} and the  Gelfand and Kolmogorov numbers by the first named author in \cite{AG}.

  The necessary and sufficient conditions for compactness of Sobolev embeddings of subspaces of Besov and Triebel-Lizorkin spaces consisted of block-radial functions were proved in \cite{LS}. The approach was based on atomic decomposition. The same tool was used to prove the Strauss inequality for block radial functions cf. \cite{STi}.  In \cite{DS} we estimates from below and from above entropy numbers of compact Sobolev embeddings of block-radial   Besov spaces  $R_\gamma B^{s}_{p,q}(\R^d)$   and fractional  Sobolev spaces $R_\gamma H^{s}_{p}(\R^d)$, $1<p<\infty$. Whereas in the present paper we study the asymptotic behaviour of  approximation numbers of these embeddings.  Our approach is based on a discretization in terms of wavelet bases of corresponding weighted function spaces and the properties of  diagonal operators acting in $l_p$ spaces. 
     
Let $SO(\gamma)= SO(\gamma_1)\times \ldots \times SO(\gamma_m)\subset SO(d)$ be a group of block-radial symmetries, $2\le \gamma_{1}\le \dots \le \gamma_{m}$,  $\gamma_1+\ldots + \gamma_m = |\gamma|= d$. We prove that the asymptotic behaviour approximation numbers  of embeddings 
\begin{equation}\label{int}
 \id: R_{\gamma} H^{s_1}_{p_1}(\Rd)\rightarrow R_{\gamma} H^{s_2}_{p_2}(\Rd) \quad \text{and}\quad \id: R_{\gamma} B^{s_1}_{p_1,q_1}(\Rd)\rightarrow R_{\gamma} B^{s_2}_{p_2,q_2}(\Rd)
 \end{equation}
  depends on the smallest dimension $\gamma_1$ of the blocks, the number $n$ of blocks with dimension $\gamma_1$ and the parameters $p_1, p_2$ but does not depends on the smoothness indices  $s_1,s_2$. More precisely the main assertion says that if $1 < p_1<p_2 < \infty$, and  $s_1-\frac{1}{p_1} -s_2 +\frac{1}{p_2}>0$, then 
	\begin{displaymath}
	a_k(\id)\,\sim\,\left\{ \begin{array}{ll} (k\log^{1-n} k)^{-\frac{\gamma_{1}-1}{p}}, &{\mathrm{if}\quad 1< p_1< p_2\leq2\ \mathrm{or}\ 2\leq p_1<p_2\leq\infty,}\\   k^{\frac{1}{t}-\frac{1}{2}}(k\log^{1-n} k)^{-\frac{\gamma_{1}-1}{p}}, &{\mathrm{if}\quad 1< p_1<2< p_2 \leq\infty \ \mathrm{and}\ \frac{\gamma_{1}}{p}>\frac{1}{\min\{p_1,p'_2\}},}\\
		(k^{t/2}\log^{1-n} k)^{-\frac{\gamma_{1}-1}{p}}, &{\mathrm{if}\quad 1< p_1<2< p_2 \leq\infty \ \mathrm{and}\ \frac{\gamma_{1}}{p}<\frac{1}{\min\{p_1,p'_2\}},} \end{array}\right.
\end{displaymath}
where $\frac{1}{t}=\frac{1}{\min\{p'_1,p_2\}}$  and $\frac{1}{p}=\frac{1}{p_1}-\frac{1}{p_2}$. 
Moreover we prove that the embedding \eqref{int} is nuclear if and only if 
$$\frac{s_1 - s_2}{d}> \frac{1}{p_1}- \frac{1}{p_2}>\frac{1}{\gamma_{1}}.$$

\subsection*{Notation}

As usual, $\N$ denotes the natural numbers, $\N_0:= \N \cup \{0\}$,  $\Z$ denotes the integers and
$\R$ the real numbers. Logarithms are always taken in base 2, $\log=\log_2$.
If $X$ and $Y$ are two (quasi)-Banach spaces, then the symbol  $X \hookrightarrow Y$ indicates that the embedding is continuous.
The set of all linear  and bounded operators
$T : X \to Y$, denoted by $\cl (X,Y)$, is  equipped with the standard norm.
As usual, the symbol  $c $ denotes positive constants
which depend only on the fixed parameters $s,p,q$ and probably on auxiliary functions,
unless otherwise stated; its value  may vary from line to line.
We will use the symbol $A \sim B$, where $A$ and $B$ can depend on certain parameters. The meaning of $A \sim B$ is given by: there exist  constants $c_1,c_2>0$ such that  inequalities $ c_1 A \le B \le c_2 A$ hold for all values of the parameters.

	If $E$ denotes a space of distributions  (functions) on $\R^d$ then by $R_\gamma E$ we mean the subset of $SO(\gamma)$-invariant  distributions (functions) in $E$. We endow this subspace  with the same norm as the original space. If $SO(\gamma)=SO(d)$, i.e. if  the subspace consists of radial functions, then we will write $RE$.  
	Similarly if $G$ is a finite group of reflections in $\R^d$ then $R_G E$  denotes  the subspace of those elements of the space $E$ that are invariant with respect to $G$. 
	 
    We assume that the reader is familiar with the basic fact concerning fractional Sobolev and Besov spaces, both unweighted and weighted with Muckenhoupt weights. We refer e.g. to  \cite{Tr1, Bui, HS} for details.  

\section{Spaces of block-radial functions}
First we recall what we mean by block-radial symmetry in $\R^d$, $d>1$. 
Let $ 1\le m\le d$ and let $\gamma\in\N^{m}$ be  an $m$-tuple  $\gamma=(\gamma_1,\ldots, \gamma_m)$, such that  $\gamma_1+\ldots + \gamma_m = |\gamma|= d$. The $m$-tuple  $\gamma$ describes the decomposition of $\R^{d}= \R^{\gamma_{1}}\times\dots\times\R^{\gamma_{m}}$  into $m$ subspaces of dimensions $\gamma_{1},\dots,\gamma_{m}$ respectively.  Let 
\[
SO(\gamma)= SO(\gamma_1)\times \ldots \times SO(\gamma_m)\subset SO(d)
\] 
be a group of isometries on $\R^{d}$. An element $g=(g_1,\ldots, g_m)$, $g_i\in SO(\gamma_i)$  acts  on $x=(\tilde{x}_1,\ldots, \tilde{x}_m)$, $\tilde{x}_i\in \R^{\gamma_i}$ by $x\mapsto g(x)= (g_1(\tilde{x}_1),\ldots, g_m(\tilde{x}_m))$.  If $m=1$ then $SO(\gamma)= SO(d)$ is the special orthogonal group acting on $\R^d$. If $m=d$ then the group  is trivial since then $\gamma_1= \ldots =\gamma_m = 1$ and $SO(1)= \{\mathrm{id}\}$.  We will always assume that $\gamma_i\ge 2$ for any $i=1,\ldots, m$. This means in particular that any group $SO(\gamma_i)$ is not trivial and  $d\ge 4$. 

Let $\bspq(\R^d)$, $s\in \R$ and $0< p,q\le \infty$, be a Besov space and $R_\gamma\bspq(\R^d)$ be its subspace consisted of $SO(\gamma)$-invariant distributions.  Similarly, we define  subspaces $R_\gamma H^s_p(\Rd)$ of fractional Sobolev spaces $H^s_p(\R^d)$, $1<p<\infty$.

It is known that  then the embeddings 
\begin{align}\label{ent2}
R_\gamma B^{s_1}_{p_1,q_1}(\R^d) \hookrightarrow R_\gamma B^{s_2}_{p_2,q_2}(\R^d) 
\quad \text{and} \quad R_\gamma H^{s_1}_{p_1}(\R^d) \hookrightarrow R_\gamma H^{s_2}_{p_2}(\R^d) 
\intertext{are compact if and only if }
s_1-s_2>d(\frac{1}{p_1}-\frac{1}{p_2})>0\qquad\text{and}\qquad \min\{\gamma_i: i=1,\ldots m\}\ge 2, \nonumber
\end{align}
 cf. \cite{LS}. 

In \cite{DS} we calculated an asymptotic behaviour of  entropy numbers of the above embeddings. In the present paper we use  essentially  the same general strategy to estimate  approximation numbers. It  consists of two steps. First  we prove using the method of traces  that the spaces $R_\gamma B^{s}_{p,q}(\R^d)$ is isomorphic to  some weighted Besov space $R_G B^{s}_{p,q}(\R^m, w_\gamma)$ defined on $\R^m$. Afterwards using the wavelet basis in the weighted space we reduce the problem to corresponding sequence spaces.   

We describe shortly these isomorphisms.  We  use  the following notation related to the action of the group $SO(\gamma)$ on $\R^d$, 
Let $\bar{\gamma}_i = 1 + \sum_{\ell=0}^{i-1} \gamma_\ell$ if $i=1,2,\ldots, m$ with ${\gamma}_0= 0$, and 
\begin{align}\nonumber
	r_{j} \,=\,r_{j}(x)\, =\, \left(x_{\bar{\gamma}_{j}}^{2}+\dots+x_{\bar{\gamma}_{j+1}-1}^{2}\right)^{1/2}, j=1,\ldots ,m. 
\end{align}
 We consider  a hyperplane 
\begin{equation}\nonumber
	H_\gamma = \mathrm{span}\{e_{\bar{\gamma}_1},e_{\bar{\gamma}_2}, , \ldots , e_{\bar{\gamma}_m}\}\ ,
\end{equation}
where $e_j$, $j=1,\ldots, d$ is a standard orthonormal basis in $\R^d$. The hyperplane $H_\gamma$ can be  identified   with $\R^m$ in the standard way so we  write $(r_1,\ldots , r_m)\in H_\gamma$ if $r_1 e_{\bar{\gamma}_1} + r_2 e_{\bar{\gamma}_2} + \ldots + r_m e_{\bar{\gamma}_m}\in H_\gamma$.  We  need also a finite group of reflections  $G=G(\gamma)$ acting on $H_\gamma$. The group consists of  transformations  $g_{i_1,\ldots,i_m}\in G(\gamma)$ given by  \[
g_{i_1,\ldots,i_m}(r_1,\ldots , r_m) = \Big( (-1)^{i_1}r_1,\dots , (-1)^{i_m}r_m\Big), \qquad  (i_1,\ldots,i_m)\in \{1,2\}^{m}\, . 
\] 
Apart of the group $G$ we consider on $H_\gamma$ also a weight $w_\gamma$  defined by the following formula   
\begin{equation} 
	w_\gamma(r_1,\dots , r_m)\  = \ \prod_{i=1}^{m}|r_i|^{\gamma_i-1}.
	\label{wg}
\end{equation}
Direct calculations show that it is a Muckenhoupt weight. The following proposition was proved in \cite{DS}   

\begin{Prop}[cf. \cite{DS}]
	\label{anfang3}
	Let $d\ge 2$, $s\in \R$, $1 < p \le \infty$ and $0<q \le \infty$.  Then the space 
	$R_{\gamma} B^s_{p,q}(\Rd)$ is isomorphic to  $R_{G} B^s_{p,q} (H_\gamma,  w_\gamma)$.  Analogously the space  
	$R_{\gamma} H^s_p (\Rd)$ is isomorphic to $R_{G}H^s_p (H_\gamma,  w_\gamma)$, if $p<\infty$. 
\end{Prop}

In consequence  we can consider the embedding
\[ R_{G} B^{s_1}_{p_1,q_1} (\R^m,  w_\gamma)\hookrightarrow R_{G} B^{s_2}_{p_2,q_2} (\R^m,  w_\gamma)\]
instead of the original one. The last embedding can be analysed via  wavelet basis of Besov spaces with Muckenhoupt weights. The existence  of the basis reduced the problem from the level of  function spaces to the level of sequence spaces  i.e. we are left with the embedding 
\[ \id:b^{\sigma_1}_{p_1,q_1}(\ w_\gamma)\rightarrow  b^{\sigma_2}_{p_2,q_2}( w_\gamma), 
 \]
where $\sigma_i=s_i-\frac{m}{p_i}+\frac{m}{2}$. 

The above sequence spaces  are defined in the following way.    Let   $Q_{\nu,k}$  be a dyadic cube in $\R^m$, centred at $2^{-\nu}k$, $k\in\Z^m$, $\nu\in\N_0$,   with the side length $2^{-\nu}$.  For $0<p<\infty$, $\nu\in\N_0$ and $k\in\Z^m$ we denote by $\chi_{\nu,k}^{(p)}$ the $p$-normalized characteristic function of the cube $Q_{\nu,k}$,
\[ 
\chi_{\nu,k}^{(p)}(x)=2^{\frac{m \nu}{p}}\chi_{\nu,k}(x)=\begin{cases}
	2^{\frac{m\nu }{p}}, &\text{for} \quad x\in Q_{\nu,k},\\  
	0, &\text{for}\quad x\notin Q_{\nu,k}.  
\end{cases}
\]
For $0<p<\infty$, $0<q\leq\infty$, $\sigma\in\R$ 
we define the   sequence spaces $b_{p,q}^\sigma(w_\gamma)$  
by
\begin{align*}
	b_{p,q}^\sigma(w_\gamma):=\biggl{\{}& \lambda=\{\lambda_{\nu,n}\}_{v,n}: \lambda_{\nu,n}\in\C,\\
	&\|\lambda|b_{p,q}^\sigma(w_\gamma)\|=
		\bigg(\sum_{\nu=0}^\infty 2^{\nu\sigma q}\bigg(\sum_{k\in\Z^m} |\lambda_{\nu,k}|^p2^{m\nu}w_\gamma(Q_{\nu,k})\bigg)^{q/p}\bigg)^{1/q} <\infty 
	\biggr{\}}.
\end{align*}
If $\sigma=0$ we write $b_{p,q}(w)$ instead of $b_{p,q}^\sigma(w)$; moreover, if $w\equiv 1$ we write $b_{p,q}^\sigma$ instead of $b_{p,q}^\sigma(w)$. One can easily prove that 
\begin{align}\label{wQ}
	w_\gamma(Q_{\nu,k}) =  \int_{Q_{\nu,k}} w_\gamma(x) dx \sim  
	2^{-\nu d} w_\gamma(Q_{0,k}) .
\end{align}
Thus to estimate the growth properties of the double-indexed sequence $\{w_\gamma(Q_{\nu,k})\}_{\nu,k}$ it is sufficient  to estimate  elements of the sequence   $\{w_\gamma(Q_{0,k})\}_{k}$. This is done in the following lemma that was proved in \cite{DS}. 
\begin{Lem}[\cite{DS}]\label{volumn}
	Let $\gamma= (\gamma_1 , \gamma_2, \ldots , \gamma_m)\in \N^m$ be a multi-index such that $\gamma_i\ge 2$ for any $i=1,\ldots ,m$.  
	We assume that   
	\begin{equation*}
		\gamma_1 \le  \gamma_2\le  \ldots \le  \gamma_m
	\end{equation*}
	and put 
	\begin{equation*}
		n = \max\{ i\,:\, \gamma_i=\gamma_1\}.
	\end{equation*} 
	
	Let  $\tau$ denote a bijection of $\Z^m$ onto $\N_0$ such that $\tau(k)< \tau(\ell)$ if $w_\gamma (Q_{0,k}) <  w_\gamma (Q_{0,\ell})$. 
	Then  there are positive constants $c_1$ and $c_2$ such that for sufficiently large $L\in \N$ the inequalities 
	\begin{equation}\label{red1}
		c_1 2^{L(d-m)} \le w_\gamma (Q_{0,k}) \,\le \,c_2 2^{L(d-m)}  
	\end{equation} 
	if and only if 
	\begin{equation*}
		c_1  2^{\frac{L(d-m)}{\gamma_1-1}}L^{n-1}\le \tau(k)  \,\le \,c_2 2^{\frac{L(d-m)}{\gamma_1-1}}L^{n-1}.
	\end{equation*}
\end{Lem}






\section{ Approximation numbers}

We start with recalling the definition of approximation numbers of bounded linear operators. 
\begin{Def}
	Let $X$ and $Y$ be quasi-Banach spaces and $T~\in~\cl(X,Y)$. 
	For $k\in \N$, we define the \textit{$k$-th approximation number} by
	\begin{displaymath}
		a_k(T):=\inf\{\|T-A\|\ :\ A\in L(X,Y),\ \rank(A)<k\},
	\end{displaymath}
	where $\rank(A)$ denotes the dimension of the range $A(X)=\{A(x): x\in X\}$.
\end{Def}
If $a_k(T)\rightarrow 0$ for $k\rightarrow \infty$, then the operator $T$ can be approximated by finite rank operators, so it is compact. The asymptotic behaviour of the sequence $a_k(T)$ gives us the quantitative information about this approximation. 
The following well known properties of approximation numbers are very useful if one calculate their  asymptotic behaviour for a given operator, cf.  \cite{ ET,Pi2,Pi}.   


\begin{Prop}\label{prop}
	Let  $W, X, Y$ be quasi-Banach spaces.
	Then
	\begin{enumerate}
		\item $ \|T\|=a_1(T)\geq a_2(T)\geq a_3(T)\geq\ldots\geq 0$\quad
		for all $T~\in~L(X,Y)$,
		\item {\normalfont{(additivity)}} $$a_{n+k-1}^p (T_1 +T_2)\leq a_k^p(T_1) + a_n^p(T_2)\quad  $$
		for all  $T_1,T_2\in L(X,Z)$ and $n,k\in\N$, if $Z$ is a $p$-Banach space. 
		\item {\normalfont(multiplicativity)} $$a_{n+k-1} (T_1T_2)\leq a_k(T_1)a_n(T_2)\quad $$
		for all $T_1\in L(X,Y)$, $T_2\in L(W,X)$ and $n,k\in\N$, 
		\item  $a_k(T)=0$ when $\rank(T)<k$
		for all $T\in L(X,Y)$. 
	\end{enumerate}
	\end{Prop}

We use  concepts of operator ideals and operator ideal quasi-norms. Both  are  very useful tools in different situations, cf. \cite{Pi, Pi2}.   Here we recall only what we need for our proofs. The most popular operators quasi-norms related to  approximation numbers are  defined via  Lorentz sequences spaces, but one can define the quasi-norm  in  more general way, cf. \cite{Kuh3, KLSS1}.    Let  $\omega=(\omega_n)$ be an increasing sequence of positive real numbers satisfying the following regularity condition $\omega_{2k}\sim\omega_k$. We define an operator quasi-norm $L_{\omega,\infty}^{(a)}$ putting   
\begin{equation*}
	L_{\omega,\infty}^{(a)}(T):=\sup _{k\in\N}\omega_k a_k(T),	
\end{equation*}
where $T\in L(X,Y)$ and $X,Y$ are quasi-Banach spaces. 
In particular if $\omega_k=k^{1/s}$ we write  
\begin{equation}\label{Ls}
	L_{s,\infty}^{(a)}(T):=\sup _{k\in\N}k^{\frac{1}{s}} a_k(T).	
\end{equation}
We will often take advantage of the fact that the quasi-norm $L_{\omega,\infty}^{(a)}$ is equivalent to an $\rho$-norm for some $\rho$, $0<\rho \leq
1$.

A class of operators for which we can calculate asymptotic behaviour of approximation numbers is a family of diagonal operators acting between $\ell_p$ spaces. 
 Let $\sigma_1\geq\sigma_2\geq...\geq 0$ be an non-increasing sequence of positive numbers and $D_\sigma:\ell_{p}\rightarrow \ell_{p}$ be a diagonal operator defined by this sequence.  It is  well-known that then 
\[a_k(D_\sigma:\ell_{p}\rightarrow \ell_{p})=\sigma_k , \qquad k\in \N. \]
cf. \cite[Proposition 2.9.5]{Pi2}
In consequence if $p_1<p_2$ then by  the monotonicity of the scale $\ell_p$-space and multiplicativity od approximation numbers  we get
\begin{equation}
	\label{Dia}
	a_k(D_\sigma:\ell_{p_1}\rightarrow\ell_{p_2})\leq\sigma_k,  \qquad k\in \N.
\end{equation}
To precisely estimate   approximation numbers of given  diagonal operators acting between different $\ell_p$ spaces,  we need exact estimates of approximation numbers of embedding of finite dimension spaces $\ell^N_p$, $N\in\N$.  We collect known results concerning these embeddings in the next lemma,   see \cite{Cae,ET,Pi2}. 

\begin{Lem}\label{Na} 
	Let $N,k\in\N$, be a positive integer and $1/t=1/\min\{p'_1,p_2\}$.
	\begin{itemize}
		\item[$(i)$] If	$	0< p_1\leq p_2\leq 2$ or 	$	2\leq p_1\leq p_2\leq\infty$, then there is a positive constant $C$ independent of $N$ and $k$ such that
		\begin{equation*}
			a_k(\id:\ell_{p_1}^N\rightarrow \ell_{p_2}^N)\leq \left\{ \begin{array}{ll} 1, & {\mathrm{if}\quad k\leq N,}
				\\
				
				0,&{\mathrm{if}\quad k>N.} \end{array}\right. 
		\end{equation*}
		
		\item[$(ii)$] If	$	0< p_1<2< p_2\leq \infty, (p_1,p_2)\neq(1,\infty)$, then there is a positive constant $C$ independent of $N$ and $k$ such that
				\begin{equation*}
			a_k(\id:\ell_{p_1}^N\rightarrow \ell_{p_2}^N)\leq C \left\{ \begin{array}{ll} 1, & {\mathrm{if}\quad k\leq N^{2/t},}
				\\
				N^{1/t}k^{-1/2} ,& {\mathrm{if}\quad N^{2/t}< k\leq N,}
				\\ 
				0,&{\mathrm{if}\quad k>N.} \end{array}\right. 
		\end{equation*}
		
		\item[$(iii)$] Let 	$	0< p_2< p_1\leq \infty$. Then
		
		\begin{equation*}
			a_k(\id:\ell_{p_1}^N\rightarrow \ell_{p_2}^N)=(N-k+1)^{\frac{1}{p_2}-\frac{1}{p_1}},\qquad k=1,...,N.
		\end{equation*}	
		
	\end{itemize}
	Moreover, if $k\leq \frac N 4$, then we have an equivalence in $(i)$ and $(ii)$.
\end{Lem}

Now we prove our main result concerning diagonal operators. 

\bp \label{mainDia}
Let $0< p_1,p_2\leq\infty$, $(p_1,p_2)\neq (1,\infty)$ and $D_\sigma$ be a diagonal operator generated by a sequence 
\[
\sigma_k=\left\{ \begin{array}{ll} k^{-\alpha}\log^\beta k, & k\ge 2^{\frac{\beta}{\alpha}}\\	
	1  &{\mathrm{otherwise},} \end{array}\right. 
\]
where $\alpha>\max(0,\frac{1}{p_2}-\frac{1}{p_1})$ and $\beta \geq 0$ . Let $\frac{1}{t}=\frac{1}{\min\{p'_1,p_2\}}$. Then  for all $k\in\N, k>1$ 
\begin{align*}
&a_k(D_\sigma :\ell_{p_1}\rightarrow \ell_{p_2})\sim \\
&\qquad \left\{ \begin{array}{ll}
	 k^{-\alpha}\log^\beta k, &{\mathrm{if}\quad 0< p_1\leq p_2\leq2\ \mathrm{or}\ 2\leq p_1\leq p_2\leq\infty,}\\  k^{-\alpha+\frac{1}{t}-\frac{1}{2}}\log^\beta k, &{\mathrm{if}\quad 0< p_1<2< p_2 \leq\infty, 
	 	 \ \mathrm{and}\ \alpha>\frac{1}{t},}\\
	k^{-\alpha\frac{t}{2}}\log^\beta k, &{\mathrm{if}\quad 0< p_1< 2< p_2 \leq\infty, 
		 \ \mathrm{and}\ \alpha<\frac{1}{t},}
	\\ k^{-\alpha+\frac{1}{p_2}-\frac{1}{p_1}}\log^\beta k, &{\mathrm{if}\quad 0< p_2< p_1 \leq\infty,} \\
\end{array}\right.
\end{align*}
\ep

\bpr
{\em Step 1. Estimation from above.} We consider different cases separately. 
The upper estimate in the case $0<p_1\leq p_2\leq2$ and $2\leq p_1\leq p_2\leq\infty$ follows from general estimates (\ref{Dia}). 
 
To consider other cases we split the diagonal operator into a sum of finite dimensional operators. 
  Let
\begin{displaymath}
	\Lambda:=\{\lambda=(\lambda_l)_l:\quad \lambda_l\in\C,\quad 1\leq l<\infty\}
\end{displaymath}
and let $P_i: \Lambda\rightarrow\Lambda$ be a projection defined in the following way
\begin{displaymath}
	P_i(\lambda)=
	\left\{ \begin{array}{ll} \lambda_l, & {\mathrm{if}\quad 2^{i-1}\leq l <2^i,} \\ 0, &{\mathrm{otherwise}}\end{array}\right.
	, \qquad l\in\N,
\end{displaymath}
for $\lambda=(\lambda_\ell)$, $\ell\in \N$,  and  $i=1,2,...$. Then
\begin{displaymath}
	D_\sigma=\sum_i D_\sigma \circ P_i.
\end{displaymath} 
Moreover for any  $M\in\N_0$ we define 
\begin{displaymath}\label{P-Q}
	P^M:=\sum_{i=1}^{M}D_\sigma\circ P_i \qquad \mathrm{and} \qquad Q_M:=\sum_{i=M+1}^{\infty}D_\sigma\circ P_i.
\end{displaymath} 
Elementary properties of the approximation numbers  yield that for any $i\in\N$
\begin{equation}
	a_k(D_\sigma \circ P_i: \ell_{p_1}\rightarrow \ell_{p_2})\leq \sigma_{2^{i-1}}a_k(id:\ell_{p_1}^{2^{i-1}}\rightarrow \ell_{p_2}^{2^{i-1}}).
	\label{eq:2}
\end{equation}

First we assume that $0< p_2< p_1 \leq\infty$.  From (\ref{eq:2}) and  Lemma \ref{Na} (iii)  we get
\begin{equation}
	L_{s,\infty}^{(a)}(D_\sigma\circ P_i: l_{p_1}\rightarrow l_{p_2})\leq c 2^{(i-1)(\frac{1}{s} + \frac{1}{p_2}-\frac{1}{p_1}-\alpha)}\max\{1,i-1 \}^\beta.
	\label{eq:9}
\end{equation} 
We choose $s>0$ such that $\frac{1}{s} + \frac{1}{p_2}-\frac{1}{p_1}-\alpha >0$. Then, from (\ref{eq:9})
\begin{equation}\label{eq:3}
	L_{s,\infty}^{(a)}(P^M)^\rho\leq \sum_{i=1}^{M}L_{s,\infty}^{(a)}(D_\sigma\circ P_i)^\rho\leq c\sum_{i=1}^{M} 2^{\rho i(\frac{1}{s} + \frac{1}{p_2}-\frac{1}{p_1}-\alpha)}i^{\rho\beta}\leq c2^{\rho M(\frac{1}{s} +\frac{1}{p_2}-\frac{1}{p_1}-\alpha)}M^{\rho\beta}.
\end{equation}
So from  (\ref{Ls}) and (\ref{eq:3}) we get 
\begin{equation}\label{eq:4}
	a_{2^M}(P^M:l_{p_1}\rightarrow l_{p_2})\leq c 2^{ M(-\alpha+ \frac{1}{p_2}-\frac{1}{p_1})}M^\beta.	
\end{equation}

On the other hand  we can choose $s>0$  such that $\frac{1}{s} + \frac{1}{p_2}-\frac{1}{p_1}-\alpha <0$ since   $\frac{1}{p_2}-\frac{1}{p_1}<\alpha$. Hence, from (\ref{eq:9})   we get 
\begin{align}\label{eq:5}
	L_{s,\infty}^{(a)}(Q_M)^\rho\leq & \sum_{i=M+1}^{\infty}L_{s,\infty}^{(a)}(D_\sigma\circ P_i)^\rho\leq c\sum_{i=M+1}^{\infty} 2^{\rho i(\frac{1}{s} + \frac{1}{p_2}-\frac{1}{p_1}-\alpha)}i^{\rho\beta}\leq\\
	 & c 2^{\rho M(\frac{1}{s} + \frac{1}{p_2}-\frac{1}{p_1}-\alpha)}M^{\rho\beta}.
	 \nonumber	
\end{align}
In view of (\ref{Ls}) and (\ref{eq:5}) we get
\begin{equation}
	a_{2^M}(Q_M:\ell_{p_1}\rightarrow \ell_{p_2})\leq c 2^{ M(-\alpha+ \frac{1}{p_2}-\frac{1}{p_1})}M^\beta.
	\label{eq:6}
\end{equation}
Therefore by (\ref{eq:4}), (\ref{eq:6}) and standard argument  
\begin{equation*}
	a_{k}(D_\sigma:\ell_{p_1}\rightarrow \ell_{p_2})\leq  c k^{ -\alpha+ \frac{1}{p_2}-\frac{1}{p_1}}\log^\beta k,
\end{equation*}
for any $k\in\N, k>1$.

Let now  $0< p_1< 2< p_2\leq\infty$, $(p_1,p_2)\neq (1,\infty)$ and $ \alpha>\frac 1 t.$ 
Then by Lemma \ref{Na} (ii) and (\ref{eq:2}) we have
\begin{equation}
	L_{s,\infty}^{(a)}(D_\sigma\circ P_i: \ell_{p_1}\rightarrow \ell_{p_2})\leq 
	c2^{(i-1)(\frac{1}{s} -\frac{1}{2}+\frac{1}{t}-\alpha)}i^\beta,\quad\mathrm{if}\ s< 2,
	\label{eq:10}
\end{equation} 
and 
\begin{equation}
	L_{s,\infty}^{(a)}(D_\sigma\circ P_i: \ell_{p_1}\rightarrow \ell_{p_2})\leq c 2^{(i-1)(\frac{1}{s}\frac{2}{t}-\alpha)}i^\beta, \quad\mathrm{if}\ s\geq 2.
	\label{eq:11}
\end{equation}
Take $s<2$ such that $\frac{1}{s} -\frac{1}{2}+\frac{1}{t}-\alpha >0$. Then the formula (\ref{eq:10}) yields
\begin{equation}
	L_{s,\infty}^{(a)}(P^M)^\rho\leq \sum_{i=1}^{M} 2^{\rho i(\frac{1}{s} -\frac{1}{2}+\frac{1}{t}-\alpha)}i^{\rho\beta}\leq c 2^{\rho M(\frac{1}{s} - \frac{1}{2}+\frac{1}{t}-\alpha)}M^{\rho\beta}.
	\label{eq:13}
\end{equation}
So we get from (\ref{Ls}) and (\ref{eq:13})
\begin{equation}\label{eq:14}
	a_{2^M}(P^M:\ell_{p_1}\rightarrow \ell_{p_2})\leq c 2^{ M(-\alpha+ \frac{1}{t}-\frac{1}{2})}M^\beta .
\end{equation}
We choose  $s=2$  and then from (\ref{eq:11}) for sufficiently large $M$
\begin{equation}	\label{eq:15}
	L_{2,\infty}^{(a)}(Q_M)^\rho\leq \sum_{i=M+1}^{\infty} 2^{\rho i(\frac{1}{t} -\alpha)}i^{\rho\beta}\leq c_2 2^{\rho M(\frac{1}{t}-\alpha)}M^{\rho\beta}.
\end{equation}
So from (\ref{Ls}) and (\ref{eq:15}) we get
\begin{equation}
	a_{2^M}(Q_M:\ell_{p_1}\rightarrow \ell_{p_2})\leq c 2^{ M(-\alpha+\frac{1}{t}-\frac{1}{2})}M^\beta.
	\label{eq:17}
\end{equation}
Hence by standard argument from (\ref{eq:14}) and (\ref{eq:17}) we have
\begin{equation*}
	a_{k}(D_\sigma:\ell_{p_1}\rightarrow \ell_{p_2})\leq c k^{-\alpha+\frac{1}{t}-\frac{1}{2}}\log^\beta k,
\end{equation*}
for any $k\in\N, k>1$.

Let now  $0< p_1< 2< p_2\leq\infty$, $(p_1,p_2)\neq (1,\infty)$ and $\alpha<\frac 1 t.$ In this case, we use a bit different  projections  $\widetilde{P_i}:\Lambda\rightarrow\Lambda$  defined   by
\begin{displaymath}
	\widetilde{P_i}(\lambda)=
	\left\{ \begin{array}{ll} \lambda_l, & {\mathrm{if}\quad [2^{\frac t 2 (i-1)}]\leq l <[2^{\frac t 2 i}],} \\ 0, &{\mathrm{otherwise}}\end{array}\right.
	, \qquad %
	\lambda=(\lambda_l), \;l\in\N, 
\end{displaymath}
where $[a]$ denotes an integer part of $a\in \R_+$, $i=1,2,...$.  Moreover we define  operators $\widetilde{P^M}$ and  $\widetilde{Q_M}$ in the analogous way to \eqref{P-Q}. 
Then the dimension of a range of the projection $\widetilde{P_i}$ is   $N_i = [2^{\frac t 2 i}]-[2^{\frac t 2 (i-1)}]\sim [2^{\frac t 2 i}]$. 
\begin{equation}\label{eq:22}
	a_k(D_\sigma \circ 	\widetilde{P_i}: \ell_{p_1}\rightarrow \ell_{p_2})\leq \sigma_{ [2^{\frac t 2 (i-1)}]}a_k(id:\ell_{p_1}^{N_i}\rightarrow \ell_{p_2}^{N_i}).
\end{equation}
So by Lemma \ref{Na} (ii) and (\ref{eq:22}) we have
\begin{equation}
	L_{s,\infty}^{(a)}(D_\sigma\circ 	\widetilde{P_i}: \ell_{p_1}\rightarrow \ell_{p_2})\leq c  2^{ i (\frac 1 s-\frac t 2\alpha)}i^\beta, \quad\mathrm{if}\ s\geq 2.
	\label{eq:31}
\end{equation}
Since $\alpha<\frac 1 t $ we can take $s\geq 2$ such that $\frac 1 s-\frac t 2\alpha>0$. Then the formula (\ref{eq:31}) yields
\begin{equation}
	L_{s,\infty}^{(a)}(\widetilde{P^M})^\rho\leq \sum_{i=1}^{M} 2^{\rho i(\frac{1}{s}-\frac{t}{2}\alpha)}i^{\rho\beta}\leq c2^{\rho M(\frac{1}{s} -\frac{t}{2}\alpha)}M^{\rho\beta}.
	\label{eq:33}
\end{equation}
So we get from (\ref{Ls}) and (\ref{eq:33})
\begin{equation}
	a_{2^M}(\widetilde{P^M}:\ell_{p_1}\rightarrow \ell_{p_2})\leq c 2^{ -M\frac t 2\alpha }M^\beta.
	\label{eq:34}
\end{equation}
When we choose $s\geq 2$ such that  $\frac{1}{s}-\frac{t}{2}\alpha <0$ then  we get from (\ref{eq:31}) that 
\begin{equation}
	L_{s,\infty}^{(a)}(\widetilde{Q_M)}^\rho\leq c_1\sum_{i=M+1}^{\infty} 2^{\rho i(\frac{1}{s}-\frac{t}{2}\alpha)}i^{\rho\beta}\leq c_2 2^{\rho M(\frac{1}{s}-\frac{t}{2}\alpha)}M^{\rho\beta}, 
	\label{eq:18}
\end{equation}
for sufficiently large $M$. Now from (\ref{Ls}) and (\ref{eq:18}) we get
\begin{equation}
	a_{2^M}(\widetilde{Q_M}:\ell_{p_1}\rightarrow \ell_{p_2})\leq c 2^{ -M\frac t 2\alpha }M^\beta.
	\label{eq:35}
\end{equation}
Hence by standard argument from (\ref{eq:34}) and (\ref{eq:35}) we have
\begin{equation*}
	a_{k}(D_\sigma:\ell_{p_1}\rightarrow \ell_{p_2})\leq c k^{ -\frac t 2\alpha }\log^\beta k,
\end{equation*}
for any integer $k, k>1$.

{\em Step 2. Estimation from below.} 
First we note that
\begin{displaymath}
	a_{k}(D_\sigma:\ell_{p_1}\rightarrow \ell_{p_2})\sim a_{k}(\id:\ell_{p_1}(\omega)\rightarrow \ell_{p_2}),
\end{displaymath}
where
\begin{displaymath}
	\ell_{p_1}(\omega)=\bigl\{(\lambda_\ell)_\ell: \left\|\lambda|\ell_{p_{1}}(\omega)\right\|=\bigl(\sum_{\ell=1}^{\infty}|
	\lambda_\ell \omega(\ell)
	|^{p_1}\bigr)^{1/p_1}<\infty\bigr\},\qquad \omega(\ell)=\sigma_\ell^{-1}.
\end{displaymath}
Thus we can consider  the following commutative diagram
\begin{equation}\nonumber
	\begin{CD}
		\ell_{p_1}^N@>{S}>>\ell_{p_1}(\omega)\\
		@V{\Id}VV@VV{\id}V\\
		\ell_{p_2}^N @<{T}<<\ell_{p_2}.\,
	\end{CD}
\end{equation}
where 
\begin{displaymath}
	(S(\lambda))_i=\left\{ \begin{array}{ll} \lambda_{i-N+1}, &{\mathrm{if}\ N\leq i\le 2N-1},\\ 0, &{\mathrm{otherwise,}}\end{array}\right.
\end{displaymath}
and
\begin{displaymath} 
	(T(\lambda))_i=\lambda_{i+N-1},\qquad\mathrm{if}\ 1\leq i\leq N.
\end{displaymath}
Then  $\|T\|=1$ and $ \|S\|= \sigma_{2N-1}^{-1}\sim \sigma_{N}^{-1} $. In consequence   multiplicativity of approximation numbers implies  
\begin{equation}\label{eq:24}
		a_k(\Id)\leq\|S\|\|T\|a_k(\id)\le  C N^\alpha\log^{-\beta}N \, a_k(\id), \qquad\mathrm{for}\; k\in\N.	
\end{equation}

In the case  $0<p_1\leq p_2\leq2$ or  $2\leq p_1\leq p_2\leq\infty$ we take $N=2^k$, $k\in\N, k\geq 2$.  Lemma \ref{Na} (i) and the inequality (\ref{eq:24})   give us 
\begin{equation}
	C\leq (\sigma_{2^k})^{-1} a_{2^{k-2}}(\id).
	\label{eq:25}
\end{equation}

Now consider the case $0<p_2<p_1\leq\infty$. We choose $N=2^{k+1}$, $k\in\N, k\geq 2$ and use  Lemma \ref{Na} (iii) and (\ref{eq:24}) to get
\begin{equation}
	C2^{(k-2)(\frac{1}{p_2}-\frac{1}{p_1})}\leq (\sigma_{2^{k-1}})^{-1} a_{2^{k-2}}(\id).
	\label{eq:28}
\end{equation}

At the end assume that  $0< p_1<2< p_2\leq\infty$, $(p_1,p_2)\neq (1,\infty)$. First let  $\frac 1 t<\alpha$. We take $N=2^k$, $k\in\N, k\geq 2$. Now by Lemma \ref{Na} (ii) and  (\ref{eq:24}) we get
\begin{equation}
	C2^{(k-2)(\frac 1 t - \frac 1 2)}\leq (\sigma_{2^k})^{-1} a_{2^{k-2}}(\id),
	\label{eq:26}
\end{equation}
for $k> \frac{2t}{t-2}$.

Let  $\frac 1 t>\alpha$. If we choose the $N=[2^{k\frac t 2}]$, $k\in\N, k\geq 2$ then
\[
N^{\frac 1 t}2^{-\frac{k-2}{2}}\sim C \quad\mathrm{and}\quad 2^{k-2}\leq \frac N 4,
\]
for sufficiently large $k$. So by  Lemma \ref{Na} (ii) and  (\ref{eq:24}) we get
\begin{equation}
	C\leq (\sigma_{[2^{k\frac{t}{2}}]})^{-1} a_{2^{k-2}}(\id).
	\label{eq:27}
\end{equation}

Now the proposition follows from  
(\ref{eq:25}), (\ref{eq:28}), (\ref{eq:26}), (\ref{eq:27}) and  the equivalence  $\sigma_k\sim\sigma_{2k}$.\epr


\begin{Prop} \label{main2seq}
	Let $\gamma= (\gamma_1 , \gamma_2, \ldots , \gamma_m)\in \N^m$, $m\in \N$,   be a multi-index such that 
	$2\le \gamma_1\le \ldots \le \gamma_m$, 
	$d= \gamma_1+ \ldots + \gamma_m$, and let $n = \max\{ i\,:\, \gamma_i=\gamma_1\}$.  
	Let $1\le p_1<p_2\le\infty$, $(p_1,p_2)\not= (1,\infty)$, $0<q_1,q_2\le \infty$ and  $s_1,s_2\in \R$.  If $\delta=s_1-s_2-d(\frac{1}{p_1}-\frac{1}{p_2})>0$,   $\frac{1}{t}=\frac{1}{\min\{p'_1,p_2\}}$  and $\frac{1}{p}=\frac{1}{p_1}-\frac{1}{p_2}$ then 
	\begin{align*}
		& a_k\big(id: b^{\sigma_1}_{p_1,q_1}( w_\gamma)\rightarrow  b^{\sigma_2}_{p_2,q_2}( w_\gamma) \big) \sim \\
		&\qquad\quad\left\{ \begin{array}{ll} (k\log^{1-n} k)^{-\frac{\gamma_{1}-1}{p}}, &{\mathrm{if}\quad 1\leq p_1< p_2\leq2\ \mathrm{or}\ 2\leq p_1<p_2\leq\infty,}\\   k^{\frac{1}{t}-\frac{1}{2}}(k\log^{1-n} k)^{-\frac{\gamma_{1}-1}{p}}, &{\mathrm{if}\quad 1\leq p_1<2< p_2 \leq\infty \ \mathrm{and}\,  
		\frac{\gamma_1}{p} > \frac{1}{\min\{p_1,p'_2\}},}\\ 
			(k^{t/2}\log^{1-n} k)^{-\frac{\gamma_{1}-1}{p}}, &{\mathrm{if}\quad 1\leq p_1<2< p_2 \leq\infty \ \mathrm{and}\ 
		\frac{\gamma_1}{p} < \frac{1}{\min\{p_1,p'_2\}},}	
	\end{array}\right.
	\end{align*}	
	where $\sigma_i= s_i +\frac{m}{2}- \frac{m}{p_i}$, $i=1,2$, $k\in\N, k>1$ .
\end{Prop}

\bpr
{\em Step 1.}
It is convenient to change slightly the notation. Let $0<p,q\le \infty$, $\sigma\in \R$ and let  ${\mathcal X}$ denote $\Z^m$ or $N_0$. We introduce the following sequence space,
\begin{align}\label{lspace}
	\ell_q(2^{\nu \sigma}\ell_p({\mathcal X}, w)) &:= \left\{\right.  \lambda=\{\lambda_{\nu,\ell}\}_{\nu,\ell}: \lambda_{\nu,\ell}\in\C, \\ 
	&\left. \|\lambda|l_q(2^{\nu\sigma}l_p({\mathcal X},w))\| = \bigg(\sum_{\nu=0}^\infty 2^{\nu\sigma q}\Big(\sum_{\ell\in {\mathcal X}} |\lambda_{\nu,\ell}|^p w_{\nu,\ell} \Big)^{q/p}\bigg)^{1/q}<\infty \right\} , \nonumber
\end{align}
where  $w= (w_{\nu,\ell})_{\nu,\ell}$, $w_{\nu,\ell} >0$. As usual we write $\ell_q(\ell_p({\mathcal X}, w))$ if $\sigma=0$ and $\ell_q(2^{\nu \sigma}\ell_p({\mathcal X}))$ if $w \equiv 1$. 

By standard arguments we get 
\begin{align} \nonumber
	a_k\big(id: b^{\sigma_1}_{p_1,q_1}(\ w_\gamma) &\rightarrow  b^{\sigma_2}_{p_2,q_2}( w_\gamma) \big) \\
	& \sim \,   
	a_k\big( id: \ell_{q_1}(\ell_{p_1}(\Z^m))\rightarrow \ell_{q_2}(2^{\nu \sigma}\ell_{p_2}(\Z^m, w^{(\gamma)}))\big), \nonumber
\end{align}
where  $\sigma = s_2- s_1$ and  $w^{(\gamma)}= \big( w^{(\gamma)}_{\nu, n}\big)$ with  $w^{(\gamma)}_{\nu, n}= w_\gamma(Q_{\nu,n})^{1-\frac{p_2}{p_1}}$.  

{\em Step 2.} 
Now we prove that  
\begin{align}\label{17.3-0}
	a_k\Big( id: \ell_{q_1}\big(\ell_{p_1}(\Z^m)\big) & \rightarrow \ell_{q_2}\big(2^{\nu \sigma}\ell_{p_2}(\Z^m, w^{(\gamma)})\big)\Big)\, \sim  \\ 
	\sim\, & a_k\Big( id: \ell_{q_1}\big(\ell_{p_1}(\N_0)\big)\rightarrow \ell_{q_2}\big(2^{-\nu \delta} \ell_{p_2}(\N_0, \widetilde{w}^{(\gamma)})\big)\Big),
	\nonumber
\end{align}
where $\widetilde{w}^{(\gamma)}_\ell = \max(1,\ell\log^{1-n)}\ell)^{(\gamma_1-1)(1-\frac{p_2}{p_1})}$.
The estimate of the weight on the cubes \eqref{wQ} and the definition of the bijection $\tau$, cf. Lemma \ref{volumn}, give us 
\begin{align} \label{17.03-1}\nonumber
	\Big\|\lambda \Big|  \ell_{q_2}&\big(2^{\nu \sigma}\ell_{p_2}(\Z^m, w^{(\gamma)})\big) \Big\|  = \\
	& = \bigg(\sum_{\nu=0}^\infty 2^{\nu\sigma q_2}\bigg(\sum_{n\in\Z^m} |\lambda_{\nu,n}|^{p_2} w_\gamma(Q_{\nu,n})^{1-\frac{p_2}{p_1}}\bigg)^{q_2/p_2}\bigg)^{1/q_2}\sim  \nonumber\\ 
	& \sim \bigg(\sum_{\nu=0}^\infty 2^{\nu(\sigma - d(\frac{1}{p_2}- \frac{1}{p_1}))q_2}\bigg(\sum_{n\in\Z^m} |\lambda_{\nu,n}|^{p_2} w_\gamma(Q_{0,n})^{1-\frac{p_2}{p_1}}\bigg)^{\frac{q_2}{p_2}}\bigg)^{\frac{1}{q_2}}\,\sim \nonumber\\
	& \sim \bigg(\sum_{\nu=0}^\infty 2^{-\nu \delta q_2}\bigg(\sum_{\ell=0}^\infty |\lambda_{\nu,\tau^{-1}(\ell)}|^{p_2} w_\gamma(Q_{0,\tau^{-1}(\ell)})^{1-\frac{p_2}{p_1}}\bigg)^{\frac{q_2}{p_2}}\bigg)^{\frac{1}{q_2}}, \nonumber
\end{align}
where $\delta = s_1-s_2 - d(\frac{1}{p_1}- \frac{1}{p_2})$. 

If $L,\ell\in \N_0$ and $2^{\frac{L(d-m)}{\gamma_1-1}}L^{n-1}\le \ell <2^{\frac{(L+1)(d-m)}{\gamma_i-1}}(L+1)^{n-1}$ then 
\begin{equation}\label{red2a}
	w_\gamma (Q_{0,\tau^{-1}(\ell)}) \,\sim\,  2^{L(d-m)}\, \sim \, (\ell\log^{1-n}\ell)^{\gamma_1-1}\, ,
\end{equation}
cf. \eqref{red1}. By what we have already proved 
\begin{align}\label{17.3-2}\nonumber
	\Big\|\lambda \Big|\ell_{q_2}\big(2^{\nu \sigma}&\ell_{p_2}(\Z^m, w^{(\gamma)})\big) \Big\|   \, \sim \, \\
	&  \sim \bigg(\sum_{\nu=0}^\infty 2^{-\nu \delta q_2}\bigg(\sum_{\ell=0}^\infty |\lambda_{\nu,\tau^{-1}(\ell)}|^{p_2} 
	\max(1,\ell\log^{1-n}\ell)^{(\gamma_1-1)(1-\frac{p_2}{p_1})}\bigg)^{\frac{q_2}{p_2}}\bigg)^{\frac{1}{q_2}}. \nonumber
\end{align}
This justifies the equivalence \eqref{17.3-0} since the estimate 
\[
\Big\|\lambda \Big|\ell_{q_1}\big(\ell_{p_1}(\Z^m)\big) \Big\|   \, \sim \, \bigg(\sum_{\nu=0}^\infty \bigg(\sum_{\ell=0}^\infty |\lambda_{\nu,\tau^{-1}(\ell)}|^{p_1} \bigg)^{\frac{q_1}{p_1}}\bigg)^{\frac{1}{q_1}}  
\]
is obvious. 

{\em Step 3.}  We prove the upper estimate of the approximation numbers. 
Let us consider the projection 
$P_\nu: \ell_{q_1}\big(\ell_{p_1}(\N_0)\big) \rightarrow \ell_{p_1}(\N_0)$
onto the $\nu$-th level,  and  the embedding operator $E_\nu: \ell_{p_2}(\N_0) \rightarrow \ell_{q_2}\big(2^{-\nu \delta} \ell_{p_2}(\N_0)\big)$, 
\[
\big(E_\nu(y)\big)_{\mu, \ell}= \begin{cases}
	y_\ell & \text{if}\qquad  \mu = \nu\\
	0 & \text{otherwise}, 
\end{cases}
\qquad y \in  \ell_{p_2}(\N_0 ).
\]
It is  obvious that 
\begin{equation}\nonumber
	\|P_\nu\| = 1 \qquad \text{and}\qquad \|E_\nu\| = 2^{-\nu\delta}\, .
\end{equation}
Let $D_\gamma$ denote the diagonal operator $D_\gamma: \ell_{p_1}(\N_0) \rightarrow \ell_{p_2}(\N_0)$ generated by the sequence
$\sigma_\ell =  (\widetilde{w}_\ell^{(\gamma)})^{\frac{1}{p_2}}$ i.e. 
$$
(D_\gamma(\lambda))_\ell = (\widetilde{w}^{(\gamma)}_\ell)^{\frac{1}{p_2}} \lambda_\ell\, .
$$ 
Then
$$
\id = \sum_{\nu=0}^\infty \id_\nu, \qquad \text{where} \qquad
\id_\nu=E_\nu D_\gamma P_\nu \,.
$$
The multiplicativity of the  approximation numbers yields
\begin{equation}
	a_k(\id_\nu) \leq c \,2^{-\nu\delta} a_k(D_\gamma)
	\label{eq:29}
\end{equation}
with a constant $c$ independent of $\nu$ and $k$. Now using  Proposition  \ref{mainDia} 
 with $\alpha=	\frac{\gamma_1-1}{p}$  we get 
\begin{align*}
&a_k(D_\gamma)\sim  \\ 
&\qquad	\left\{ \begin{array}{ll} (k\log^{1-n} k)^{-\frac{\gamma_{1}-1}{p}}, &{\mathrm{if}\quad 1\leq p_1< p_2\leq2\ \mathrm{or}\ 2\leq p_1< p_2\leq\infty,}\\  k^{\frac 1 t -\frac 1 2}(k\log^{1-n} k)^{-\frac{\gamma_{1}-1}{p}}, &{\mathrm{if}\quad 1\leq p_1<2< p_2 \leq\infty, 
 \ \mathrm{and}\ 	\frac{\gamma_1}{p} > \frac{1}{\min\{p_1,p'_2\}},}\\ (k^{t/2}\log^{1-n}k)^{-\frac{\gamma_{1}-1}{p}}, & {\mathrm{if}\quad 1\leq p_1<2< p_2 \leq\infty, 
  \ \mathrm{and}\ 	\frac{\gamma_1}{p} < \frac{1}{\min\{p_1,p'_2\}}.}		
	\end{array}\right.
\end{align*}
Now using (\ref{eq:29}) and taking $\omega_k=(a_k(D_\gamma))^{-1}$ we have
$$
L_{\omega,\infty}^{(a)}(\id_\nu)=\sup _k \, \omega_k a_k(\id_\nu) \leq c
\,2^{-\nu\delta}  
$$
and then

$$
L_{\omega,\infty}^{(a)}(\id)^r \leq  \sum_{\nu=0}^\infty
L_{\omega,\infty}^{(a)}(\id_\nu)^r \leq c \sum_{\nu=0}^\infty 2^{-\nu \delta r} < \infty \, , 
$$
which proves the upper estimate.

{\em Step 4.}  In this step we  estimate the approximation numbers from below. We consider  the following commutative diagram
\begin{equation}\nonumber
	\begin{CD}
		\ell_{p_1}^{N}@>{T}>>\ell_{q_1}\big(\ell_{p_1}(\N_0)\big)\\
		@V{\Id}VV@VV{\id}V\\
		\ell_{p_2}^{N} @<{S}<<\ell_{q_2}\big(2^{-\nu\delta}\ell_{p_2}(\N_0, \tilde{w}^{(\gamma)})\big).\,
	\end{CD}
\end{equation}
Here the operators $S$ and $T$ are defined by  
\begin{align*}
	\big(T(\xi_1, ... , \xi_N)\big)_{\nu,\ell} = 
	\begin{cases} 
		\xi_{\ell + 1-N}  & \text{if}\;  \nu=0 \;\text{and}\; N \leq l \leq 2N-1,\\
		0  &\text{otherwise} 
	\end{cases}
\end{align*}
and
$$
S((\lambda_{\nu,\ell})_{\nu,\ell} ) = (\lambda_{0,N}, \cdots ,\lambda_{0,2N-1}). 
$$

The norms of the above operators have the obvious estimates 
$$
\|T\| \leq 1    \quad {\rm and } \quad   \|S\| \leq (\widetilde{w}^{(\gamma)}_N)^{-\frac{1}{p_2}} \, .
$$
Using the multiplicativity of the approximation numbers
we get,   
\begin{align*}
	&   a_k(\Id: \ell_{p_1}^N \rightarrow \ell_{p_2}^N) \nonumber \\[0.1cm]
	& \leq \|S\|  \,a_k(\id: \ell_{q_1}(\ell_{p_1}(\N_0)) \rightarrow      \nonumber
	\ell_{q_2}(2^{-\nu \delta}\ell_{p_2}(\N_0, \widetilde{w}^{(\gamma)}))) \,\|T\| \le \\
	& (N\log^{1-n} N)^{(\gamma_1-1)(\frac{1}{p_1}-\frac{1}{p_2})}\,a_k(\id: \ell_{q_1}(\ell_{p_1}(\N_0)) \rightarrow      \nonumber
	\ell_{q_2}(2^{-\nu \delta}\ell_{p_2}(\N_0, \widetilde{w}^{(\gamma)}))).\,  
\end{align*}
We get similarly inequality as (\ref{eq:24}) in the proof Proposition \ref{mainDia}  (where now $\alpha=(\gamma_1-1)(\frac{1}{p_1}-\frac{1}{p_2})$, $\beta=(1-n)(\gamma_1-1)(\frac{1}{p_1}-\frac{1}{p_2})$). So when we take the same $N$ which was need to get (\ref{eq:25}), (\ref{eq:28}), (\ref{eq:26}), (\ref{eq:27}) in respective cases we prove the proposition.
\epr

From the above proposition  and  arguments used  in \cite{DS} for entropy numbers,  we finally arrive at the following theorem.

\begin{T}\label{mainT2}
	Let $\gamma= (\gamma_1 , \gamma_2, \ldots , \gamma_m)\in \N^m$, $m\in \N$,   be a multi-index such that $2\le \gamma_1\le \ldots \le \gamma_m$, 
	$d= \gamma_1+ \ldots + \gamma_m$, and let $n = \max\{ i\,:\, \gamma_i=\gamma_1\}$.	
	
	Let $1 < p_1<p_2 \leq \infty$,   $0<q_1,q_2\le \infty$ and  $s_1,s_2\in \R$.  If  $\delta=s_1-s_2-d(\frac{1}{p_1}-\frac{1}{p_2})>0$, $\frac{1}{t}=\frac{1}{\min\{p'_1,p_2\}}$  and $\frac{1}{p}=\frac{1}{p_1}-\frac{1}{p_2}$ then 
	\begin{align*}
	& a_k\Big( \id: R_{\gamma} B^{s_1}_{p_1,q_1}(\Rd)\rightarrow R_{\gamma} B^{s_2}_{p_2,q_2}(\Rd)\Big) \sim  \\  
	&\qquad\qquad 	\left\{ \begin{array}{ll} (k\log^{1-n} k)^{-\frac{\gamma_{1}-1}{p}}, &{\mathrm{if}\quad 1< p_1< p_2\leq2\ \mathrm{or}\ 2\leq p_1<p_2\leq\infty,}\\   k^{\frac{1}{t}-\frac{1}{2}}(k\log^{1-n} k)^{-\frac{\gamma_{1}-1}{p}}, &{\mathrm{if}\quad 1< p_1<2< p_2 \leq\infty \ \mathrm{and}\ 	\frac{\gamma_1}{p} > \frac{1}{\min\{p_1,p'_2\}},}\\
			(k^{\frac{t}{2}}\log^{1-n} k)^{-\frac{\gamma_{1}-1}{p}}, &{\mathrm{if}\quad 1< p_1<2< p_2 \leq\infty \ \mathrm{and}\ 	\frac{\gamma_1}{p} < \frac{1}{\min\{p_1,p'_2\}},} \end{array}\right.
		\end{align*}
		for $k\in\N, k>1$.
		
	   If $1 < p_1<p_2 < \infty$  then the same estimates hold for $	a_k\Big( \id: R_{\gamma} H^{s_1}_{p_1}(\Rd)\rightarrow R_{\gamma} H^{s_2}_{p_2}(\Rd)\Big) $. 	
\end{T}
\begin{Rem} {\rm 
	If $m=1$, $\gamma_1=d$ then we regain the estimate of approximation numbers of embedding of subspaces of radial functions, 
	\begin{align*}
	a_k \sim 	\left\{ \begin{array}{ll} 
	k^{-\frac{d-1}{p}}, &{\mathrm{if}\quad 1< p_1< p_2\leq2\ \mathrm{or}\ 2\leq p_1<p_2\leq\infty,}\\   
	k^{\frac{1}{t}-\frac{1}{2}-\frac{d-1}{p}}, &{\mathrm{if}\quad 1< p_1<2< p_2 \leq\infty \ \mathrm{and}\	\frac{d}{p} > \frac{1}{\min\{p_1,p'_2\}},}\\
	k^{-\frac{t(d-1)}{2p}}, &{\mathrm{if}\quad 1< p_1<2< p_2 \leq\infty \ \mathrm{and}\ 	\frac{d}{p} < \frac{1}{\min\{p_1,p'_2\}},} 
	\end{array}\right.
	\end{align*}
	cf. \cite{ST}. The proof of Proposition 1 in \cite{ST} contains  a gap in the case  $p_1<2< p_2 \leq\infty$ and  $\frac{d}{p} < \frac{1}{\min\{p_1,p'_2\}}$. This is corrected in our proof of Proposition 3.   }
\end{Rem}




\section{Nuclear embeddings}

\subsection{Nuclear operators} 
Now we prove the sufficient and necessary condition for nuclearity of the embeddings we work with. Let us recall the definition and basic facts concerning nuclear operators.  

\begin{Def}
Let $X,Y$ be Banach spaces, $T\in\cl (X,Y)$. Then $T$ is called nuclear, denoted by  $T\in\cn (X,Y)$, if there exist elements $a_j\in X^*$, the dual space of $X$, and $y_j\in Y, j\in \N$, such that $\sum_{j=0}^{\infty}\|a_j\|_{X^*}\|y_j\|_Y<\infty$ and a nuclear representation $Tx=\sum_{j=0}^{\infty}a_j(x)y_j$ for any $x\in X$. Together with the nuclear norm
$$v(T)=\inf \big\{\sum_{j=0}^{\infty}\|a_j\|_{X^*}\|y_j\|_Y:T=\sum_{j=0}^{\infty}a_j(\cdot)y_j\big\},$$

where the infimum is taken over all nuclear representations of $T$, the space $\cn (X,Y)$ becomes a Banach space.
\end{Def}

The nuclear operators have been introduced by A.Grothendieck \cite{Gr}
and were studied by many authors afterwards, cf. \cite{Pi2,Pi3,Pi} and also \cite{Pi4} for some history. It is obvious that any nuclear operator can be approximated by finite rank operators, hence they are, in particular, compact.

The following properties of nuclear operators are  needed later on.
\begin{Prop}
	\label{Pro1}
\begin{itemize}
	\item[(i)]
If $X$ is an $n$-dimensional Banach space, then
$$
v (id:X\rightarrow X)=n.
$$
\item[(ii)]
For any Banach space $X$ and any bounded linear operator $T:l_\infty^n \rightarrow X$ we have
$$ v(T)=\sum_{i=1}^n\|Te_i\|.$$
\item[(iii)]
If 
$T\in \cl(X,Y)$ is a nuclear operator and $S\in \cl(X_0,X)$ and $R\in \cl(Y,Y_0)$, then $STR$ is a nuclear operator and
$$
v(STR)\leq \|S\|\|R\|v(T).
$$
\end{itemize}
\end{Prop}

Once more we will need estimates for diagonal operators. We introduce the following notation: for numbers $r_1,r_2\in \left[ 1,\infty\right]$, let $t(r_1,r_2)$ be given by
\begin{equation*}
\frac{1}{t(r_1,r_2)}=\begin{cases}
1, &\text{if} \quad 1\leq r_2\leq r_1\leq \infty,\\  
1-\frac{1}{r_1}+\frac{1}{r_2}, &\text{if}\quad  1\leq r_1\leq r_2\leq \infty.  
\end{cases}
\end{equation*}
Hence $1\leq t(r_1,r_2) \leq\infty$, and
$$
\frac{1}{t(r_1,r_2)}=1-\big(\frac{1}{r_1}-\frac{1}{r_2}\big)_+\geq \frac{1}{r^*}=\big(\frac{1}{r_2}-\frac{1}{r_1}\big)_+,
$$
with $t(r_1,r_2)=r^*$ if, and only if, $\{r_1,r_2\}=\{1,\infty\}$.\\
Recall that $c_0$ denotes the subspace of $l_\infty$ of all sequences converging to null. The following proposition was proved by A.Tong. 

\begin{Prop}[{\cite[Thms. 4.3, 4.4]{To}}]
	\label{Pro2}	
	 Let $1\leq r_1, r_2 \leq\infty$ and $D_\tau$ be the above diagonal operator.
	 \begin{itemize}
	 	\item[(i)] Then $D_\tau$ is nuclear if, and only if, $\tau=(\tau_j)_j\in l_{t(r_1,r_2)}$, with $l_{t(r_1,r_2)}=c_0$ if $t(r_1,r_2)=\infty$. Moreover,
	 	$$
	 	v(D_\tau:l_{r_1}\rightarrow l_{r_2})=\|\tau|l_{t(r_1,r_2)}\|.
	 	$$
\item[(ii)] Let $n\in\N$ and $(D_\tau^n:l_{r_1}^n\rightarrow l_{r_2}^n)$ be the corresponding diagonal operator $D_\tau^n: x=(x_j)^n_{j=1}\mapsto (\tau_j x_j)_{j=1}^n$. Then
$$
v(D_\tau^n:l_{r_1}^n\rightarrow l_{r_2}^n)=\|\tau|l_{t(r_1,r_2)}^n\|.
$$

	 	\end{itemize}
\end{Prop}

\subsection{Nuclear embeddings of spaces of block-radial functions}

In this section we find sufficient and necessary condition for nuclearity of the embeddings \eqref{int}. 

\begin{Prop} \label{main1seq}
	Let $\gamma= (\gamma_1 , \gamma_2, \ldots , \gamma_m)\in \N^m$, $m\in \N$,   be a multi-index such that 
	$2\le \gamma_1\le \ldots \le \gamma_m$, 
	$d= \gamma_1+ \ldots + \gamma_m$, and let $n = \max\{ i\,:\, \gamma_i=\gamma_1\}$.  
	Let $1\le p_1,p_2\le\infty$, $0<q_1,q_2\le \infty$ and  $s_1,s_2\in \R$.  Then the embedding 
	\begin{equation*}
id: b^{\sigma_1}_{p_1,q_1}( w_\gamma)\rightarrow  b^{\sigma_2}_{p_2,q_2}( w_\gamma),
	\end{equation*}
	where $\sigma_i= s_i +\frac{m}{2}- \frac{m}{p_i}$, $i=1,2$ is nuclear if, and only if
	\begin{equation}
	\label{c1}
 \frac{s_1 - s_2}{d}> \frac{1}{p_1}- \frac{1}{p_2}>\frac{1}{\gamma_{1}}. 
	\end{equation}
\end{Prop}

\bpr
{\em Step 1.}
We use the   sequence space introduced in Step 1 and Step 2 of the proof of Proposition \ref{main2seq} cf. \eqref{lspace}-\eqref{red2a}. By the same argument as above we can prove that  
\begin{align*}
&id: b^{\sigma_1}_{p_1,q_1}(\ w_\gamma) \rightarrow  b^{\sigma_2}_{p_2,q_2}( w_\gamma)  \intertext{is nuclear if  and  only  if}  
&id: \ell_{q_1}\big(\ell_{p_1}(\N_0)\big)\rightarrow \ell_{q_2}\big(2^{-\nu \delta} \ell_{p_2}(\N_0, \widetilde{w}^{(\gamma)})
\end{align*}
is nuclear, where $\widetilde{w}^{(\gamma)}_\ell = \max(1,\ell\log^{1-n}\ell)^{(\gamma_1-1)(1-\frac{p_2}{p_1})}$ and  $\delta = s_1-s_2 - d(\frac{1}{p_1}- \frac{1}{p_2})$.  

{\em Step 2.}  We first deal with sufficiency of  (\ref{c1}) for the nuclearity.
Let us consider the projection 
$P_\nu: \ell_{q_1}\big(\ell_{p_1}(\N_0)\big) \rightarrow \ell_{p_1}(\N_0)$
onto the $\nu$-th vector-coordinate,  and  the embedding operator $E_\nu: \ell_{p_2}(\N_0) \rightarrow \ell_{q_2}\big(2^{-\nu \delta} \ell_{p_2}(\N_0)\big)$, 
\[
\big(E_\nu(y)\big)_{\mu, \ell}= \begin{cases}
y_\ell & \text{if}\qquad  \mu = \nu\\
0 & \text{otherwise}, 
\end{cases}
\qquad y \in  \ell_{p_2}(\N_0).
\]
It is  obvious that 
\begin{equation}\nonumber
\|P_\nu\| = 1 \qquad \text{and}\qquad \|E_\nu\| = 2^{-\nu\delta}\, .
\end{equation}
Let $D_\gamma$ denote the diagonal operator $D_\gamma: \ell_{p_1}(\N_0) \rightarrow \ell_{p_2}(\N_0)$ generated by the sequence
$\sigma_\ell =  (\widetilde{w}_\ell^{(\gamma)})^{\frac{1}{p_2}}$ i.e. 
\begin{equation}
(D_\gamma(\lambda))_\ell = (\widetilde{w}^{(\gamma)}_\ell)^{\frac{1}{p_2}} \lambda_\ell\, .
\label{D}
\end{equation}
Then
\begin{equation}
\id = \sum_{\nu=0}^\infty \id_\nu, \qquad \text{where} \qquad
\id_\nu=E_\nu D_\gamma P_\nu \,.
\label{v1.1}
\end{equation}
The multiplicativity of the nuclear norm ( Proposition \ref{Pro1} (iii)) yields
\begin{equation}
v(\id_\nu) \leq c \,2^{-\nu\delta} v(D_\gamma).
\label{v1.2}
\end{equation}
Now using  Proposition \ref{Pro2} (i)
 \begin{equation}
v(D_\gamma)= \|\sigma_l |\ell_{t(p_1,p_2)}\|.
\label{v1.3}
\end{equation}

So it only remains to calculate the last norm. First assume that $t(p_1,p_2)<\infty$. In this case,

\begin{align}
\nonumber
 \|\sigma_l |\ell_{t(p_1,p_2)}\|^{t(p_1,p_2)}=\sum_{l=0}^\infty \max(1,\ell\log^{1-n}\ell)^{(\gamma_1-1)(\frac{1}{p_2}-\frac{1}{p_1})t(p_1,p_2)}\leq \\ c \sum_{l=2}^\infty (\ell\log^{1-n}\ell)^{(\gamma_1-1)(\frac{1}{p_2}-\frac{1}{p_1})t(p_1,p_2)}<\infty
 \label{v1.4}
\end{align}

since $(\gamma_1-1)(\frac{1}{p_1}-\frac{1}{p_2})t(p_1,p_2)>1$, cf. (\ref{c1}). Now by (\ref{v1.1}), (\ref{v1.2}), (\ref{v1.3}), (\ref{v1.4}) we get
$$
v(\id) \leq c\sum_{\nu=0}^\infty \,2^{-\nu\delta},
$$
where $c>0$ is independent of $\nu$. By (\ref{c1}) we have that  $\delta>0$ therefore 
$v(id)<\infty$.

If $t(p_1,p_2)=\infty$, i.e., if $p_1=1$ and $p_2=\infty$, then

 $$v(D_\gamma)=\sup_{\ell\geq 0}\max(1,\ell\log^{1-n}\ell)^{(1-\gamma_1)}<\infty$$
since  $\gamma_1>1$  and in consequence similarly as in the previous case $v(id)<\infty$.

{\em Step 3.} 
Now we show the necessity of (\ref{c1}) for the nuclearity of $id$. Since any nuclear map is compact, the nuclearity of $id$ implies its compactness which by (\ref{ent2}) leads to 
$ \frac{s_1 - s_2}{d}> \frac{1}{p_1}- \frac{1}{p_2}>0$ (so $p_1<p_2$).
Now let us assume that $\frac{1}{p_1}- \frac{1}{p_2}\leq  \frac{1}{\gamma_{1}}$ and  consider the following commutative diagram
\begin{equation}\nonumber
\begin{CD}
\ell_{p_1}(\N_0)@>{T}>>\ell_{q_1}\big(\ell_{p_1}(\N_0)\big)\\
@V{\Id}VV@VV{\id}V\\
\ell_{p_2}(\N_0,\tilde{w}^{(\gamma)}) @<{S}<<\ell_{q_2}\big(2^{-\nu\delta}\ell_{p_2}(\N_0, \tilde{w}^{(\gamma)})\big).\,
\end{CD}
\end{equation}
Here the operators $S$ and $T$ are defined by  
\begin{align*}
\big(T(\xi_j)\big)_{\nu,\ell} = 
\begin{cases} 
\xi_j  & \text{if}\;  \nu=0 \;\text{and}\;  l \in\N_0 ,\\
0  &\text{otherwise} 
\end{cases}
\end{align*}
and
$$
S((\lambda_{\nu,\ell})_{\nu,\ell} ) = (\lambda_{0,l})_l. 
$$
The multiplicativity of the nuclear norm ( Proposition \ref{Pro1} (iii)) yields
$$
   v(\Id) \leq \|S\|v(\id) \|T\| \le v(\id),
$$
since  
$$
\|T\| \leq 1    \quad {\rm and } \quad   \|S\| \leq 1 \, .
$$
But
$$
v(\Id: \ell_{p_1}(\N_0)\rightarrow \ell_{p_2}(\N_0,\tilde{w}^{(\gamma)}))\sim  v(D_\gamma:\ell_{p_1}(\N_0)\rightarrow \ell_{p_2}(\N_0)),
$$
where $ D_\gamma$ denote the diagonal operator defined in (\ref{D}) i.e.
$$
(D_\gamma(\lambda))_\ell = \max(1,\ell\log^{1-n}\ell)^{(\gamma_1 - 1)(1/p_2-1/p_1)} \lambda_\ell\, 
$$ 
and Proposition \ref{Pro2} (i) implies
$$v(D_\gamma)= \|\sigma_l |\ell_{t(p_1,p_2)}\|.$$
So if $t(p_1,p_2)<\infty$ then 
\begin{equation}
\label{nr1}
\bigg(\sum_{l=0}^\infty \max(1,\ell\log^{1-n}\ell)^{(\gamma_1-1)(\frac{1}{p_2}-\frac{1}{p_1})t(p_1,p_2)}\bigg)^{\frac{1}{t(p_1,p_2)}}= v(D_\gamma)\leq v(\id),
\end{equation}
for arbitrary $\ell\in\N$ and
 if $t(p_1,p_2)=\infty$ ($p_1=1, p_2=\infty$)
\begin{equation}
\label{nr2}
(\ell\log^{1-n}\ell)^{(1-\gamma_1)}\leq\sup_{\ell\geq 0}\max(1,\ell\log^{1-n}\ell)^{(1-\gamma_1)}= v(D_\gamma)\le  v (\id),
\end{equation}
for arbitrary $\ell\in\N$. 

But $(\gamma_1-1)(\frac{1}{p_1}-\frac{1}{p_2})t(p_1,p_2)\leq 1$ since $0<\frac{1}{p_1}- \frac{1}{p_2}\leq\frac{1}{\gamma_1}$. So the series (\ref{nr1}) is divergent. Similarly $1-\gamma_1\geq 0$ if  $t(p_1,p_2)=\infty$, so the supremum in (\ref{nr2})  is infinite. But this lead to a contradiction to $ v(\id)<\infty$. 
Thus  $\frac{1}{p_1}- \frac{1}{p_2}> \frac{1}{\gamma_1}$.
 This proves the proposition. \epr

\begin{Prop} \label{main1function}
	Let $\gamma= (\gamma_1 , \gamma_2, \ldots , \gamma_m)\in \N^m$, $m\in \N$,   $2\le \gamma_1\le \ldots \le \gamma_m$, 
	$d= \gamma_1+ \ldots + \gamma_m$, and let $n = \max\{ i\,:\, \gamma_i=\gamma_1\}$. 	
	Let $1\le p_1,p_2\le\infty$, $0<q_1,q_2\le \infty$ and  $s_1,s_2\in \R$. Then the embedding 
	\begin{equation*}
	\id :R_G B^{s_1}_{p_1,q_1}(\R^m, w_\gamma) \hookrightarrow R_G B^{s_2}_{p_2,q_2}(\R^m, w_\gamma)
	\end{equation*}
	is nuclear if, and only if
	\be
	\label{nuBG}
 \frac{s_1 - s_2}{d}> \frac{1}{p_1}- \frac{1}{p_2}>\frac{1}{\gamma_{1}}. 
\ee
	
\end{Prop}

\bpr
{\em Step 1.  Sufficient condition.} 
Proposition \ref{main1seq} and the wavelet decomposition of the spaces give us 
	\begin{equation*}
\Id :B^{s_1}_{p_1,q_1}(\R^m, w_\gamma)\rightarrow B^{s_2}_{p_2,q_2}(\R^m, w_\gamma)
\end{equation*}
is nuclear if, and only if
\be
\label{nuB}
 \frac{s_1 - s_2}{d}> \frac{1}{p_1}- \frac{1}{p_2}>\frac{1}{\gamma_{1}}. 
\ee

So it remains to show that we have the same sufficient condition for the $G(\gamma)$-invariant subspaces. The operator 
$$
P f (x)= \frac{1}{|G(\gamma)|}\sum_{g\in G(\gamma)} f(g(x)) 
$$
is a bounded projection of $B^{s}_{p,q}(\R^m, w_\gamma)$ onto $R_G B^{s}_{p,q}(\R^m, w_\gamma)$ since $G(\gamma)$ is a finite group of linear  isometries.  So 
the  sufficient condition \eqref{nuBG} follows from \eqref{nuB} and  the following  commutative diagram 
\begin{equation}\nonumber
\begin{CD}
R_GB^{s_1}_{p_1,q_1}(\R^m, w_\gamma) @ >{\id}>> R_GB^{s_2}_{p_2,q_2}(\R^m, w_\gamma)\\
@V{id}VV                                          @AA{P}A \\
B^{s_1}_{p_1,q_1}(\R^m, w_\gamma) @ >>{\Id}> B^{s_2}_{p_2,q_2}(\R^m, w_\gamma). 
\end{CD}
\end{equation}

{\em Step 2.  Necessary condition.} 
Since any nuclear map is compact, the nuclearity of $\id$ implies its compactness which by \eqref{ent2} leads to $s_1 - s_2> d(\frac{1}{p_1}- \frac{1}{p_2})>0$ (so $p_1<p_2$).
Now assume that $\gamma_1(\frac{1}{p_1}- \frac{1}{p_2})\leq1$.
To prove the necessary condition we use also the wavelet decomposition. The group $G(\gamma)$ divides $\R^m$ into the finite sum of cones, that have pairwise disjoint interiors. Let us choose one of those cones and denote it by $\widetilde{\mathcal{C}}$. 
Moreover, let $\eta$ be a function belonging to $C^\ell(\Rd)$, $\ell>s_1$, such that $\supp \eta \subset \widetilde{\mathcal{C}}$ and $\eta(x)=1$ if $x\in {\mathcal{C}} = \{x\in \widetilde{\mathcal{C}}: \dist (x, \partial\widetilde{\mathcal{C}})>\varepsilon\}$, for some fixed sufficiently small $\varepsilon>0$. 
We consider a family  $\mathcal{K} = \{k\in \Z^m:\; \supp \phi_{0,k} \subset \mathcal{C}\}$. Here  $\phi_{0,k}$ are elements of the wavelet basis generated by the father wavelet (scaling function), cf. \cite[Appendix]{DS}.  The group $G(\gamma)$ is finite therefore for any $0<c_1<c_2$ we can find   $L_0$ such that for any $L\ge L_0$, $L\in \N$, we have 
\begin{align}\label{below1}
\#\{k\in \mathcal{K}: & \;c_1 2^{L(d-m)}\le w(Q_{0,k})\le c_2 2^{L(d-m)} \} \, \sim\,  \\ 
\sim\, &  \#\{k\in \Z^d:   c_1 2^{L(d-m)}\le w(Q_{0,k})\le c_2 2^{L(d-m)} \} .  \nonumber
\end{align}   
Please note that $L$ large means that the cubes $Q_{0,k}$ are located far from the origin, cf. \eqref{wQ}. The above considerations and Lemma \ref{volumn} yield  the existence of the bijection $\sigma: \mathcal{K} \rightarrow \N_0$ such that 
\begin{equation}\label{red1bis}\nonumber
w_\gamma (Q_{0,k}) \,\sim\,  2^{L(d-m)} \qquad \Longleftrightarrow\qquad  \sigma(k)\, \sim \, 2^{\frac{L(d-m)}{\gamma_1-1}}L^{n-1} \, .
\end{equation}

Let $v_{\gamma}(\ell) =  w_\gamma (Q_{0,\sigma^{-1}(\ell)})$.  Then \eqref{red2a} and \eqref{below1} imply
\begin{equation}\label{below1a}
v_\gamma(\ell)  
\, \sim \, (\ell\log^{1-n}\ell)^{\gamma_1-1}\, .
\end{equation}  

For  further arguments we need three linear bounded  operators: $T:\ell_{p_1}(\N_0,v_\gamma)\rightarrow \ell_{q_1}(\ell_{p_1}(\Z^m, w_\gamma))$, $M_\eta: B^{s_2}_{p_2,q_2}(\R^m, w_\gamma)\rightarrow B^{s_2}_{p_2,q_2}(\R^m, w_\gamma)$ and $S:\ell_{q_2}(\ell_{p_2}(\Z^m, w_\gamma))\rightarrow \ell_{p_2}(\N_0,v_\gamma)$. The operators are defined in the following way:
\begin{align*}
(T\lambda)_{j,k} &\, = \begin{cases}
\lambda_{\sigma(k)}& \quad\text{if}\quad k\in \mathcal{K}\quad \text{and}\quad j=0,\\
0 & \quad \text{otherwise},
\end{cases}
\\
M_\eta(f) &\, = |G(\gamma)| \eta \cdot f, \qquad f\in  B^{s_2}_{p_2,q_2}(\R^m, w_\gamma) ; \\
S(\lambda)_\ell &\, = \lambda_{0, \sigma^{-1}(\ell)} \qquad \text{it}\qquad \lambda\in  \ell_{q_2}(\ell_{p_2}(\Z^m, w_\gamma))\, .
\end{align*}

Using these operators we can construct the following commutative diagram 
\begin{equation}\label{diag230317}\nonumber 
\begin{CD}
\ell_{p_1}(\N_0,v_\gamma)@>{T}>> \ell_{q_1}(\ell_{p_1}(\Z^m, w_\gamma))     @>{\mathcal{W}^{-1}}>> B^{s_1}_{p_1,q_1}(\R^m, w_\gamma) @ >{P}>>  R_GB^{s_1}_{p_1,q_1}(\R^m, w_\gamma) \\ 
@V{\Id}VV     &{} &{} &{}  &{}                                                        @VV{\id}V\\
\ell_{p_2}(\N_0,v_\gamma)@<S<<  \ell_{q_2}(\ell_{p_2}(\Z^m, w_\gamma)) @<{\mathcal{W}}<< B^{s_2}_{p_2,q_2}(\R^m, w_\gamma) @ <{M}<<  R_GB^{s_2}_{p_2,q_2}(\R^m, w_\gamma).
\end{CD}
\end{equation}

Here   $\mathcal{W}$ is the isomorphism defined by the wavelet basis.  It follows from the above diagram that  
\begin{equation}\label{below2}
v(\Id: \ell_{p_1}(\N_0,v_\gamma)\rightarrow \ell_{p_2}(\N_0, v_\gamma))\le C  v (\id).
\end{equation} 
But 
\begin{equation}\label{below3}
v(\Id: \ell_{p_1}(\N_0,v_\gamma)\rightarrow \ell_{p_2}(\N_0,v_\gamma))\, \sim \, 
v(D_\gamma: \ell_{p_1}(\N_0)\rightarrow \ell_{p_2}(\N_0)),
\end{equation} 
where  $D_\gamma$ denotes the diagonal operator  generated by the sequence
$\sigma_\ell =  (\ell\log^{1-n}\ell)^{(\gamma_1 - 1)(1/p_2-1/p_1)}$ i.e. 
$$
(D_\gamma(\lambda))_\ell = (\ell\log^{1-n}\ell)^{(\gamma_1 - 1)(1/p_2-1/p_1)} \lambda_\ell\, .
$$ 
Using once more  Proposition \ref{Pro2} (i) and \eqref{below1a}-\eqref{below3} we have
$$   \big(\sum_{l=2}^\infty (\ell\log^{1-n}\ell)^{(\gamma_1-1)(\frac{1}{p_2}-\frac{1}{p_1})t(p_1,p_2)}\big)^{1/t(p_1,p_2)} \sim v(D_\gamma)\le  v (\id),$$
if $t(p_1,p_2)<\infty$ and
$$\sup_{\ell\geq 2}(\ell\log^{1-n}\ell)^{(1-\gamma_1)}\sim v(D_\gamma)\le  v (\id),$$
if $t(p_1,p_2)=\infty$ ($p_1=1, p_2=\infty$).

However, the above series is divergent to infinity if $0<\frac{1}{p_1}-\frac{1}{p_2}\leq \frac{1}{\gamma_{1}}$. 
 Also, the sequence $(\ell\log^{1-n}\ell)^{(1-\gamma_1)}$ is divergent to infinity if $\gamma_{1}\leq 1$.
But this lead to contradiction for $ v(\id)<\infty$.
 Thus $\frac{1}{p_1}-\frac{1}{p_2}> \frac{1}{\gamma_{1}}$ 
which ends the proof.

\epr


Summarizing the previously given arguments, we have the following main theorem.
\begin{T} \label{main2function}
	Let $\gamma= (\gamma_1 , \gamma_2, \ldots , \gamma_m)\in \N^m$, $m\in \N$,   $2\le \gamma_1\le \ldots \le \gamma_m$, 
	$d= \gamma_1+ \ldots + \gamma_m$, and let $n = \max\{ i\,:\, \gamma_i=\gamma_1\}$. 	
	Let $1< p_1,p_2\le\infty$, ($1<p_1, p_2<\infty$ in the case of Sobolev spaces)  $0<q_1,q_2\le \infty$ and  $s_1,s_2\in \R$. Then the embedding 
	\begin{equation*}
R_\gamma B^{s_1}_{p_1,q_1}(\R^d) \hookrightarrow R_\gamma B^{s_2}_{p_2,q_2}(\R^d)
\qquad \bigg(R_\gamma H^{s_1}_{p_1}(\R^d) \hookrightarrow R_\gamma H^{s_2}_{p_2}(\R^d)\bigg)
	\end{equation*}
	is nuclear if, and only if
	$$
\frac{	s_1-s_2}{d}>\frac{1}{p_1}-\frac{1}{p_2}> \frac{1}{ \gamma_{1}}. 
	$$
\end{T}
\begin{Rem}
	\rm{If $m=1$ and $d=\gamma_1$ then we  get the sufficient and necessary conditions for nuclearity of embbedings of radial functions proved in \cite{HS2}
		$$
	\frac{	s_1-s_2}{d}>\frac{1}{p_1}-\frac{1}{p_2}> \frac{1}{d}. 
	$$
 }
\end{Rem}



		

		





\bigskip\bigskip~
\small

\noindent Alicja Dota\\
\noindent Institute of Mathematics, Pozna\'n University of Technology, \\
\noindent ul. Piotrowo 3A,  60-965 Pozna\'n\\ 
\noindent {\it E-mail}:  \texttt{alicja.dota@put.poznan.pl}

\bigskip

\noindent Leszek Skrzypczak\\
\noindent  Faculty of Mathematics and Computer Science,
Adam Mickiewicz University, \\
\noindent ul. Uniwersytetu Pozna\'nskiego 4, 61-614 Pozna\'n,
Poland\\
\noindent {\it E-mail}:  \texttt{lskrzyp@amu.edu.pl}

\end{document}